\documentclass[a4paper,12pt,reqno]{amsart}

\usepackage{amsmath,amssymb,amsfonts,latexsym}

\newtheorem{theorem}{Theorem}[section]
\newtheorem{definition}{Definition}[section]
\numberwithin{equation}{section}

\usepackage[margin=1.5 in]{geometry} 
\usepackage[colorlinks=true,
linkcolor=blue,     
citecolor=red,      
urlcolor=magenta,   
filecolor=cyan]{hyperref} 

\begin{document}
	\begin{center}{\bf \large
		{ON RECENT PARTITION FUNCTION OF \\[0.15cm] KAUR AND RANA}  }\vspace{0.2cm}
	\end{center}

	\begin{center}
		\footnotemark
		\bf Anjelin Mariya Johnson and \footnotetext{Email: anjelinvallialil@pondiuni.ac.in}
		\footnotemark
		\bf S. N. Fathima  \footnotetext{Email: dr.fathima.sn@pondiuni.ac.in}

	\end{center}
	{\baselineskip .5cm \begin{center}Ramanujan School of Mathematical Sciences,\\Department of Mathematics,\\ Pondicherry University,\\ Puducherry - 605 014, India.\\
	\end{center}}
	
	\noindent{\bf Abstract:} 
		Recently, Kaur and Rana introduced the partition function denoted by $\rho(n)$, where the largest part $\lambda$ appears exactly once, and the remaining parts constitute a partition of $\lambda$. In this paper, we establish new generating functions for certain variants of $\rho(n)$. Further, we obtain a linear recurrence relation for our new generating function. \\
	\vspace{.125cm}
	
	\vspace{.05cm}


	
	\section {Introduction}
	Throughout this paper, we adopt the standard notations on partitions and $q$-series, as in Andrews \cite{r1} and Gasper and Rahman \cite{r5} respectively. The $q$- shifted factorial $(a;q)_n$ is defined by 
	\begin{align*}
		(a;q)_n =
		\begin{cases}
			1&,  \;\text{for} \;\;n = 0 \\
			\prod\limits_{k=0}^{n-1} \left(1-aq^k\right)&,  \; \text{for } \;n \geq1,
		\end{cases}
	\end{align*}
	where
	$(a;q)_\infty= \lim\limits_{n \to \infty} (a;q)_n=\prod\limits_{k=0}^{\infty} \left(1-aq^k\right).$\\
	Since the infinite product diverges when $a\neq0$ and $|q|\geq1$, whenever $(a;q)_\infty$ appears in an identity, we shall assume $|q|<1$.\\
	
	Recall that a partition of a positive integer $n$ is a non-increasing sequence of positive integers $\lambda_1,\lambda_2,\ldots \lambda_n$, whose sum is $n$. Each $\lambda_i$ is called a part of the partition.
Let $p(n)$ denote the number of partitions of $n$ (see \cite{r18}, A000041]). The generating function for $p(n)$ is given by 
\begin{align*}
	\sum_{n=0}^{\infty}p(n)q^n= \frac{1}{(q;q)_\infty },
\end{align*}
with the usual convention that $p(0)=1$. Several prominent mathematicians have contributed to the study of partitions. For a general overview of theory of partitions, we refer the reader to the monumental book of Andrews  \cite{r1}.

\vspace{1em}
By imposing certain restrictions on the parts of the partition, one can obtain variants of the partition function.
For example, a partition of $n$ is {\it $\ell$-regular} if none of its parts are multiples of $\ell$. Let $b_\ell(n)$ denote the number of $\ell$-regular partitions of $n$. The 3-regular partitions of 5 are\\
                \[5\;,\;4+1\;,\;2+2+1\;,\;2+1+1+1\;,\;1+1+1+1+1.\] \\
Using elementary techniques, the generating function for $b_\ell(n)$ is given by (see \cite{r14})
\begin{align*}
	\sum_{n=0}^{\infty}b_\ell(n)q^n&=\frac{(q^\ell;q^\ell)_\infty}{(q;q)_\infty}.
\end{align*}
Interestingly, in classical representation theory the number of irreducible p-modular representations of the symmetric group $S_n$ is same as $b_p(n)$, where $p$ is prime (see \cite{r12},\cite{r8}).\\

In \cite{r4}, Corteel and Lovejoy introduced the overpartition function $\overline{p}(n)$, which counts the number of partitions of $n$ wherein the first occurrence of parts may be overlined. For example, $\overline{p}(4)=14$, since the partitions in question are\\ $4\;,\;\overline{4}\;,\;3+1\;,\;3+\overline{1}\;,\;\overline{3}+1\;,\;\overline{3}+\overline{1}\;,\;2+2\;,\;\overline{2}+2\;,\;2+1+1\;,\;\overline{2}+1+1\;,\;2+\overline{1}+1\;,\;\overline{2}+\overline{1}+1\;,\;1+1+1+1\;,\;\overline{1}+1+1+1.$\\
The generating function for $\overline{p}(n)$ is given by\\
\begin{align*}
	\sum_{n=0}^{\infty}\overline{p}(n) q^{n}\;&= \frac{(q^{2};q^{2})_{\infty}}{(q;q)_{\infty}^{2}}.
\end{align*}

Further, Lovejoy \cite{r13} investigated the $\ell$-regular overpartition $\overline{b_\ell}(n)$, which counts the number of overpartitions of $n$ with no parts divisible by $\ell$. From the above example, it is clear that $\overline{b}_3(4)=10$.
The generating function for $\overline{b_\ell}(n)$ is given by
\begin{align*}
	\sum_{n=0}^{\infty} \overline{b}_{\ell}(n) q^{n}\;\;&=\;\; \frac{(q^{\ell};q^{\ell})^2_{\infty}\, (q^{2};q^{2})_{\infty}}{\, (q;q)_{\infty}^{2}\, (q^{2\ell};q^{2\ell})_{\infty}}.\\
\end{align*}
Similarly, the number of overpartitions of $n$ in which only odd parts are used is denoted by $\overline{po}(n)$, and the number of overpartitions of $n$ in which only even parts are used is denoted by $\overline{pe}(n)$. Hence $\overline{po}(4)=6$ and $\overline{pe}(4)=4$. The generating functions for $\overline{po}(n)$ and $\overline{pe}(n)$ are given by(see \cite{r16},\cite{r6})
\begin{align*}
	\sum_{n=0}^{\infty} \overline{po}(n) q^{n}\;&= \frac{(q^{2};q^{2})_{\infty}^{3}}{(q;q)_{\infty}^{2} (q^{4};q^{4})_{\infty}},\\
	and\;\;\;
	\sum_{n=0}^{\infty}\overline{pe}(n)q^{n}\;&=\;\;\;\frac{(q^{4};q^{4})_{\infty}}{(q^{2};q^{2})_{\infty}^{2}},
\end{align*}
respectively. A part in a partition is said to have $k$ distinct colors if each part in the partition is allowed with $k$ different copies (see\cite{r10}). Let $p_{-k}(n)$ denote the number of partitions of $n$ with each parts having $k$ different colors. The generating function for
$p_{-k}(n)$ is 
\begin{align*}
	\sum_{n=0}^{\infty} p_{-k}(n) q^n &= \frac{1}{(q^k;q^k)_\infty}.
\end{align*}
For instance, if each part of partition of 3 have colors, say red(r) and blue(b) then $p_2(3)=10$, with the corresponding partitions

\;\;\;\;\;$\;\;\;\;3_r, \;3_b, \;2_r+1_r,\; 2_r+1_b,\; 2_b+1_r, \;2_b+1_b,\\1_r+1_r+1_r,\;1_r+1_r+1_b,\; 1_r+1_b+1_b,\; 1_b+1_b+1_b$\\

In \cite{r3} Chan investigated cubic partition $a(n)$, which counts the number of partition  in which the even parts can occur in two distinct colors. The generating function for $a(n)$ is given by
\begin{align*}
 	 \sum_{n=0}^{\infty} a(n) q^n\;\;\;&=\frac{1}{(q;q)_{\infty}\, (q^{2};q^{2})_{\infty}}.
 \end{align*}

	Recently, Hirschhorn and sellers \cite{r7} studied the POD function, which counts the number of partitions of $n$ wherein the odd parts are distinct (and the even parts are unrestricted). The generating function for $pod(n)$ is
	\begin{align*}
		\sum_{n=0}^{\infty} pod(n)\, q^{n}&=\frac{(q^{2};q^{2})_{\infty}}{(q;q)_{\infty}\, (q^{4};q^{4})_{\infty}}.
	\end{align*}
	Further, Andrews, Hirschhorn and sellers \cite{r2} studied the PED function, which counts the number of partitions of $n$ wherein the even parts are distinct (and the odd parts are unrestricted). The generating function for $ped(n)$ is
	\begin{align*}
		\sum_{n=0}^{\infty} ped(n)\, q^{n}&=\frac{(q^{4};q^{4})_{\infty}}{(q;q)_{\infty}}.
	\end{align*}
	\vspace{1em} 
	
	Very recently, Kaur and Rana \cite{r9} introduced the partition function $\rho(n)$ where the largest part appears exactly once, and the remaining parts constitute a partition of that largest part. For example, 
	$\rho(12)=10$, and the relevant partitions are 
	
	$6+5+1\;,\;6+4+2 \;,\; 6+4+1+1\; ,\; 6+3+3\; ,\; 6+3+2+1 \;,\; 6+3+1+1+1,\\6+2+2+2\; , \;6+2+2+1+1\; ,\; 6+2+1+1+1+1 \;,\; 6+1+1+1+1+1+1.$\\
	
	The generating function for the partition $\rho(n)$ is given by 
	\begin{align}
		\sum_{n=0}^{\infty}\rho(n)q^n&=\frac{1}{(q^2,q^2)_\infty}-\frac{1}{(1-q^2)}.\label{eq 1.1}
	\end{align}
	In this paper, motivated by the results of Kaur and Rana, we aim to investigate generating functions  for different variants of $\rho(n)$. To state our main results, we consider the following partition functions.
	\begin{definition}\label{d:2.1}
		For a positive integer $n$, we define the partition function 
		\begin{itemize}
			\item $\rho_\ell(n)$, which counts the number of partitions of $n$, wherein the largest part $\lambda$ appears exactly once, and the remaining parts constitute  $\ell$ -regular partitions of $\lambda$ .
			\item $\overline{\rho}(n)$, which counts the number of partitions of $n$, wherein the largest part $\lambda$ appears exactly once, and the remaining parts constitute overpartitions of $\lambda$.
			\item $\overline{\rho_o}(n)$, which counts the number of partitions of $n$, wherein the largest part $\lambda$ appears exactly once, and the remaining parts constitute overpartitions of $\lambda$ into odd parts.
			\item $\overline{\rho_e}(n)$, which counts the number of partitions of $n$, wherein the largest part $\lambda$ appears exactly once, and the remaining parts constitute overpartitions of $\lambda$ into even parts.
			\item $\overline{\rho}_\ell(n)$, which counts the number of partitions of $n$, wherein the largest part $\lambda$ appears exactly once, and the remaining parts constitute $\ell$-regular overpartitions of  $\lambda$ .
			\item $\rho_{pod}(n)$, which counts the number of partitions of $n$, wherein the largest part $\lambda$ appears exactly once, and the remaining parts constitute distinct, odd partitions of $\lambda$.
			\item $\rho_{ped}(n)$, which counts the number of partitions of $n$, wherein the largest part $\lambda$ appears exactly once, and the remaining parts constitute distinct, even partitions of $\lambda$.
			\item $\rho_{-k}(n)$, which counts the number of partitions of $n$, wherein the largest part $\lambda$ appears exactly once, and the remaining parts constitute $k$-coloured partition of $\lambda$.
			\item $\overline{\rho}_c(n)$, which counts the number of partitions of $n$, wherein the largest part $\lambda$ appears exactly once, and the remaining parts constitute cubic partitions of $\lambda$.

		\end{itemize}
		
	\end{definition}
	We now present our main results
	
	
	\begin{theorem}\label{t:2.1.1}
	For $n \geq 0$,
		\begin{align*}
			\sum_{n=0}^{\infty}\rho_\ell(n)q^n&= \frac{(q^{2\ell};q^{2\ell})_\infty}{(q^2;q^2)_\infty}-\frac{1}{1-q^2}+\frac{q^{2\ell}}{1-q^{2\ell}}.
			\end{align*}
	\end{theorem}

		
	\begin{theorem}\label{t:2.1.2}
	 For $n \geq 0$,  
	\begin{align}
		\sum_{n=0}^{\infty}\overline{\rho}(n)q^n&= \frac{(q^{4};q^{4})_\infty}{(q^2;q^2)^2_\infty}-\frac{2}{1-q^2}+1,\label{eq 1.2}\\
		\sum_{n=0}^{\infty}\overline{\rho}_o(n)q^n&= \frac{(q^{4};q^{4})^3_\infty}{(q^2;q^2)^2_\infty{(q^8;q^8)_\infty}}-\frac{2q^2}{1-q^4}-1,\label{eq 1.3}\\
		\sum_{n=0}^{\infty}\overline{\rho}_e(n)q^n&= \frac{(q^{8};q^{8})_\infty}{(q^4;q^4)^2_\infty}-\frac{2q^4}{1-q^4}-1.\label{eq 1.4}
	\end{align}
	\end{theorem}
	\begin{theorem}\label{t:2.1.3}
	For $n \geq 0$, 
	\begin{align*}
			\sum_{n=0}^{\infty}\overline{\rho_\ell}(n)q^n&= \frac{{(q^{2\ell};q^{2\ell})_\infty^2}{(q^4;q^4)}_\infty}{{(q^2;q^2)_\infty^2}{(q^{4\ell};q^{4\ell})_\infty}}-\frac{2}{1-q^2}+\frac{2q^{2\ell}}{1-q^{2\ell}}+1.
	\end{align*}
	\end{theorem}
	\begin{theorem}\label{t:2.1.4}
For $n \geq 0$, 
	\begin{align*}
		\sum_{n=0}^{\infty}{\rho_{-k}}(n)q^n&=\frac{1}{(q^2;q^2)^k}-\frac{k{q^2}}{1-q^2}-1.
		\end{align*}
\end{theorem}
	
	\begin{theorem}\label{t:2.1.5}
			For $n \geq 0$,  
		\begin{align*}
			\sum_{n=0}^{\infty}\rho_{c}(n)q^n&=\frac{1}{(q^2;q^2)_\infty(q^4;q^4)_\infty}-\frac{2}{1-q^2}+\frac{q^6}{1-q^4}+1+q^2.
		\end{align*}
	\end{theorem}
		\begin{theorem}\label{t:2.1.6}
		For $n \geq 0$,  
		\begin{align}
			\sum_{n=0}^{\infty}\rho_{pod}(n)q^n&=\frac{(q^4;q^4)_\infty}{{(q^2;q^2)_\infty}{(q^8;q^8)_\infty}}-\frac{1}{1-q^2},\label{eq 1.5}\\
			\sum_{n=0}^{\infty}\rho_{ped}(n)q^n&=\frac{(q^8;q^8)_\infty}{(q^2;q^2)_\infty}-\frac{1}{1-q^2}.\label{eq 1.6}
		\end{align}
	\end{theorem}
	\begin{theorem}\label{t:2.1.7}
		For $n \geq 0$, let $\rho_{\epsilon}(n)$ enumerates the number of partitions of $n$ in which the largest part say $r$ appears exactly once and the remaining parts are partitions of r where every even part is less than each odd part, (For example, $\rho_{\epsilon}(10)=3$ with the relevant partitions being 
		$5+3+2,\;5+3+1+1,\;5+1+1+1+1+1$), then
	\begin{align*}
		\sum_{n=0}^{\infty}\rho_{\epsilon}(n)q^n&=\frac{1}{1-q^2}\left[\frac{1}{(q^4;q^4)_\infty}-1\right].
		\end{align*}
	\end{theorem}

	Our proofs presented in section 2 are elementary in nature relying on generating function manipulations. We conclude this paper with an interesting recurrence relation involving partition function $\rho(n)$.
	\section{Proof of Theorems }
	
	\begin{proof}[\it\textbf{Proof of Theorem \ref{t:2.1.1}}]
		We have 
		 \begin{align*}
			\sum_{n=0}^{\infty}\rho_\ell(n)q^n&= \sum_{\substack{n=2 \\ \ell | n}}^{\infty}q^n(b_\ell(n))q^n+\sum_{\substack{n=2 \\ \ell \nmid n}}^{\infty}q^n(b_\ell(n)-1)q^n\\
			&=\sum_{\substack{n=2 \\ \ell | n}}^{\infty}(b_\ell(n))q^{2n}+\sum_{\substack{n=2 \\ \ell \nmid n}}^{\infty}(b_\ell(n)-1)q^{2n}\\
			&=\sum_{n=0}^{\infty}(b_\ell(n)-1)q^{2n}+\sum_{\substack{n=2 \\ \ell | n}}^{\infty}q^{2n}\\
			&=\sum_{n=0}^{\infty}(b_\ell(n))q^{2n}-\frac{1}{1-q^2}+\sum_{k=1}^{\infty}q^{2\ell k}.
		\end{align*}
		This completes the proof of Theorem  1.1.
	\end{proof}
		\begin{proof}[\it\textbf{Proof of Theorem \ref{t:2.1.2}}]
		We have 
			\begin{align*}
					\sum_{n=0}^{\infty}\overline{\rho}(n)q^n&=q^2(q^{1+1}+q^{\overline{1}+1})+q^3(q^{2+1}+q^{\overline{2}+1}+q^{2+\overline{1}}+q^{\overline{2}+\overline{1}})+\ldots\\
					&=\sum_{n=2}^{\infty}q^n(\overline{p}(n)-2)q^n\\
					&=1+\sum_{n=0}^{\infty}\overline{p}(n)q^{2n}-\frac{2}{1-q^2}.
				\end{align*}
				This completes the proof of (1.2).
				 The proof of (1.3) and (1.4) is similar to (1.2), Hence we omit.
			\end{proof}
		\begin{proof}[\it\textbf{Proof of Theorem \ref{t:2.1.3}}]
			We have
	\begin{align*}
		\sum_{n=0}^{\infty}\overline{\rho_\ell}(n)q^n&=\sum_{\substack{n=2 \\ \ell | n}}^{\infty}q^n\overline{b}_\ell(n)q^n+\sum_{\substack{n=2 \\ \ell \nmid n}}^{\infty}q^n(\overline{b}_\ell(n)-2)q^n\\
		&=\sum_{n=0}^{\infty}(b_\ell(n)-2)q^{2n}+2\sum_{\substack{n=2 \\ \ell | n}}^{\infty}q^{2n}+1\\
		&=\sum_{n=0}^{\infty}(b_\ell(n))q^{2n}-\frac{2}{1-q^2}+\sum_{k=1}^{\infty}q^{2\ell k}+1.\\
	\end{align*}
	This completes the proof of Theorem 1.3.		
\end{proof}
		\begin{proof}[\it\textbf{Proof of Theorem \ref{t:2.1.4}}]	
		We have
			\begin{align*}
				\sum_{n=0}^{\infty}{\rho_{-k}}(n)q^n&=\sum_{n=2}^{\infty}q^n({p_{-k}}(n)-k)q^n\\
				&=\sum_{n=2}^{\infty}({p_{-k}}(n)-k)q^{2n}\\
				&=\sum_{n=0}^{\infty}{p_{-k}}(n)q^{2n}-\frac{k}{1-q^2}+{q^2}k-1.
				\end{align*}	
		This completes the proof of Theorem 1.4.
		\end{proof}
	\begin{proof}[\it\textbf{Proof of Theorem \ref{t:2.1.5}}]	
	We have
	\begin{align*}
		\sum_{n=0}^{\infty}\rho_{c}(n)q^n&=q^2(q^{1+1})+q^3(q^{{2_1}+1}+q^{{2_2}+2}+q^{1+1+1})+\ldots\\
		&=\sum_{ \substack {n=2 \\ 2|n}}^{\infty}q^n(a(n)-2)q^n+\sum_{ \substack {n=3 \\ 2\nmid n}}^{\infty}q^n(a(n)-1)q^n\\
		&=\sum_{n=2}^{\infty}(a(n)-2)q^{2n}+\sum_{ \substack {n=3 \\ 2\nmid n}}^{\infty}q^{2n}\\
		&=\sum_{n=0}^{\infty}a(n)q^{2n}-\frac{2}{1-q^2}+\frac{q^6}{1-q^4}+1+q^2.
	\end{align*}	
	This completes the proof of Theorem 1.5.
\end{proof}
\begin{proof}[\it\textbf{Proof of Theorem \ref{t:2.1.6}}]	
		\begin{align*}
			\sum_{n=0}^{\infty}\rho_{pod}(n)q^n&=q^3(q^{2+1})+q^4(q^{3+1}+q^{2+2})+q^5(q^{4+1}+q^{3+2}+q^{2+2+1})+\ldots\\
			&=\sum_{n=3}^{\infty}q^n(pod(n)-1)q^n\\
			&=\sum_{n=0}^{\infty}{pod(n)}q^{2n}-\sum_{n=3}^{\infty}q^{2n}.
		\end{align*}
		This completes the proof of (1.5), and the proof of (1.6) is similar to (1.5). We omit the remaining  proof of Theorem 1.7.
	\end{proof}

		\section{Recurrence relation involving $\rho(n)$ and $\rho_a(n)$}
In this section we recall the counting function $a(n)$ studied by Merca \cite{r15}. For a positive integer $n$, $a(n)$ is defined to be the sum of parts counted without multiplicity in all the partitions of $n$. For example, $a(3)=7$.\\

To obtain our recurrence relation involving $\rho(n)$ and $\rho_a(n)$  we highlight the generating function satisfied by a(n).
\begin{theorem}\label{3.1} {\textnormal{(\cite{r15},Theorem 1.2)}} 
We have \\
\begin{align*}
	\sum_{n=1}^{\infty} a(n)q^n&= \frac{1}{(q;q)_\infty}.\frac{q}{(1-q)^2}.
\end{align*}

\end{theorem}
\begin{theorem}\label{3.2}
\begin{align*}
	2\rho_a(n)=n(\rho(n)-1)+2a(n/2).
\end{align*}

\end{theorem}
Proof: Using \eqref{eq 1.1} and Theorem \ref{3.1} we complete the proof of Theorem 3.2. 


\end{document}